# q-Fibonacci polynomials and q-Genocchi numbers

Johann Cigler

Fakultät für Mathematik, Universität Wien
johann.cigler@univie.ac.at

**Abstract**
We show that Genocchi and Bernoulli numbers are closely related to Fibonacci polynomials and derive some $q-$analogues.

## 1. Fibonacci polynomials and Genocchi numbers

Define the sequence $(G_{2n})_{n\geq 1} = (1,1,3,17,155,2073,38227,929569,\cdots)$ of Genocchi numbers $G_{2n}$ by their exponential generating function

$$\frac{2z}{1+e^z} = z + z\frac{1-e^z}{1+e^z} = \sum_{n\geq 0} g_n \frac{z^n}{n!} = z + \sum_{n\geq 1}(-1)^n G_{2n}\frac{z^{2n}}{(2n)!}. \tag{1.1}$$

It is well-known that $G_{2n} = (-1)^n 2(1-2^{2n})B_{2n}$, where $(B_n)$ is the sequence of Bernoulli numbers defined by

$$B_n = \sum_{k=0}^{n}\binom{n}{k}B_k \tag{1.2}$$

for $n > 1$ with $B_0 = 1$.

Let

$$F_n(s) = \sum_{k=0}^{n-1}\binom{n-k-1}{k}s^k \tag{1.3}$$

denote the Fibonacci polynomials. They are characterized by the recursion

$$F_n(s) = F_{n-1}(s) + sF_{n-2}(s) \tag{1.4}$$

with initial values $F_0(s) = 0$ and $F_1(s) = 1$ and are explicitly given by the Binet formula

$$F_n(s) = \frac{\alpha^n - \beta^n}{\alpha - \beta} \tag{1.5}$$

with $\alpha = \dfrac{1+\sqrt{1+4s}}{2}$ and $\beta = \dfrac{1-\sqrt{1+4s}}{2} = 1-\alpha.$

They satisfy the identity

$$\sum_{n\geq 0}\frac{F_n(s)}{n!}z^n = -e^z\sum_{n\geq 0}\frac{F_n(s)}{n!}(-z)^n. \tag{1.6}$$

For (1.6) is equivalent with $e^{\alpha z} - e^{\beta z} = -e^z \left( e^{-\alpha z} - e^{-\beta z} \right)$, which is trivially true.

An easy consequence is

$$(1+e^z) \sum_{n \geq 0} \frac{F_{2n}(s)}{(2n)!} z^{2n} = (e^z - 1) \sum_{n \geq 0} \frac{F_{2n+1}(s)}{(2n+1)!} z^{2n+1}. \tag{1.7}$$

If we define the linear functional $L$ on $\mathbb{C}[s]$ by

$$L\left(F_{2k+1}(s)\right) = [k=0], \tag{1.8}$$

we get

**Theorem 1.1 ( D. Dumont and J. Zeng [4], Corollary 1)**

$$L(F_{2n}(s)) = (-1)^{n-1} G_{2n}. \tag{1.9}$$

In order to show this apply $L$ to (1.7). This gives

$$\sum_{n \geq 0} \frac{L(F_{2n}(s))}{(2n)!} z^{2n} = -z \frac{1-e^z}{1+e^z} = z - \frac{2z}{1+e^z} = \sum_{n \geq 1} (-1)^{n-1} G_{2n} \frac{z^{2n}}{(2n)!}.$$

**Corollary 1.2**

$$g_n = -L(F_n(s)) \tag{1.10}$$

*for $n \neq 1$.*

As another consequence we get

**Corollary 1.3**

$$F_{2n}(s) = \sum_{k=0}^{n-1} a(n,k) F_{2k+1} \tag{1.11}$$

*with*

$$a(n,k) = (-1)^{n-k-1} \frac{1}{2k+1} \binom{2n}{2k} G_{2n-2k}. \tag{1.12}$$

E.g. we have

$$\left(a(n,k)\right)_{1\leq n\leq 5, 0\leq k\leq 4} = \begin{pmatrix} 1 & 0 & 0 & 0 & 0 \\ -1 & 2 & 0 & 0 & 0 \\ 3 & -5 & 3 & 0 & 0 \\ -17 & 28 & -14 & 4 & 0 \\ 155 & -255 & 126 & -30 & 5 \end{pmatrix}.$$

The proof follows by writing (1.6) in the form

$$\sum_{n\geq 0} \frac{F_{2n}(s)}{(2n)!} z^{2n} = -z \frac{1-e^z}{1+e^z} \sum_{n\geq 0} \frac{F_{2n+1}(s)}{(2n+1)!} z^{2n} \tag{1.13}$$

and comparing coefficients.

## 2. Fibonacci polynomials and Bernoulli numbers

In an analogous way we define a linear functional $M$ by

$$M(F_{2n}) = [n = 1]. \tag{2.1}$$

Then we get

$$M(F_{2n+1}(s)) = (2n+1)B_{2n}. \tag{2.2}$$

This follows from (1.7) and the well-known identity

$$\frac{z^2}{2} \frac{e^z + 1}{e^z - 1} = \sum_{n\geq 0} \frac{B_{2n}}{(2n)!} z^{2n+1}.$$

This implies as above

$$F_{2n+1}(s) = \sum_{j=0}^{n} \binom{2n+1}{2j+1} \frac{B_{2n-2j}}{j+1} F_{2j+2}(s). \tag{2.3}$$

The first terms of the sequence $\left((2n+1)B_{2n}\right)_{n\geq 0}$ are
$1, \frac{1}{2}, -\frac{1}{6}, \frac{1}{6}, -\frac{3}{10}, \cdots.$

We can now give a simple proof of

**Theorem 2.1 (A. v. Ettingshausen [5], L. Seidel [9], M. Kaneko [7])**

$$\sum_{i=0}^{n+1} \binom{n+1}{i} (n+i+1) B_{n+i} = 0. \tag{2.4}$$

**Proof**
We need the identity
$$\sum_{k=0}^{n}(-1)^{n-k}\binom{n}{k}F_{n+k}(s) = \sum_{k=0}^{n}(-1)^{k}\binom{n}{k}F_{2n-k}(s) = 0 \tag{2.5}$$

which is immediate from $\sum_{k=0}^{n}(-1)^{k}\binom{n}{k}\left(\alpha^{2n-k} - \beta^{2n-k}\right) = \left(\alpha^2 - \alpha\right)^n - \left(\beta^2 - \beta\right)^n = s^n - s^n = 0.$

This implies
$$\sum_{i=0}^{n+1}\binom{n+1}{2i}F_{2n+2-2i}(s) = \sum_{i=0}^{n+1}\binom{n+1}{2i+1}F_{2n+1-2i}(s).$$

If we apply the linear functional $M$ to (2.5) we get

$$\sum_{i=0}^{n+1}\binom{n+1}{2i+1}(2n-2i+1)B_{2n-2i} = 0 \text{ for } n>1.$$ Since $B_{2i+1} = 0$ for $i>0$ we have also

$$\sum_{i=0}^{n+1}\binom{n+1}{i}(n+i+1)B_{n+i} = 0 \text{ for } n>1.$$ It is easy to verify that (2.4) holds for $n=0$ and $n=1$ too.

**Remark 2.2**
This theorem has first been proved by A. v. Ettingshausen [5] and has been rediscovered by L. Seidel [9], VIII, and by M. Kaneko [7]. Therefore I will call it **v. Ettingshausen-Seidel-Kaneko identity**.

The proof by v. Ettingshausen starts with the definition of the Bernoulli numbers
$\sum_{k\geq 0}\binom{r}{k}B_k = B_r$ for $r \geq 2$ and $B_0 = 1$. He computes the differences $\Delta^w B_r$ for both sides using

the fact that $\Delta^w\binom{r}{k} = \binom{r}{k-w}$. This gives

$$\sum_{i=0}^{w}(-1)^{w-i}\binom{w}{i}B_{r+i} = \Delta^w B_r = \sum_{k=0}^{r}B_k\Delta^w\binom{r}{k} = \sum_{k=0}^{r}B_k\binom{r}{k-w} = \sum_{j=0}^{r}\binom{r}{j}B_{w+j}.$$

Choosing $r = w = n$ gives $\sum_{i=0}^{n}(-1)^{n-i}\binom{n}{i}B_{n+i} = \sum_{i=0}^{n}\binom{n}{i}B_{n+i}$ or equivalently

$$\sum_{i}\binom{n}{2i+1}B_{2n-2i-1} = 0, \tag{2.6}$$

which implies the well-known result $B_{2i+1} = 0$ for $i \geq 1$.
For $w = n, r = n+1$ he gets

$$\sum_{i=0}^{n}(-1)^{n-i}\binom{n}{i}B_{n+1+i} = \sum_{j=0}^{n+1}\binom{n+1}{j}B_{n+j} = B_n + \sum_{i=0}^{n}\binom{n+1}{i+1}B_{n+i+1}. \tag{2.7}$$

Since $B_{n+1+i} = 0$ for $n+i \equiv 0 \pmod 2$ this identity is the same as

$$-\sum_{i=0}^{n}\binom{n}{i} B_{n+1+i} = B_n + \sum_{i=0}^{n}\binom{n+1}{i+1} B_{n+i+1}. \text{ Because of } \binom{n+1}{i+1}+\binom{n}{i} = \frac{n+i+2}{n+1}\binom{n+1}{i+1}$$

this gives $(n+1)B_n + \sum_{i=0}^{n}\binom{n+1}{i+1}(n+i+2)B_{n+i+1} = 0$ or equivalently

$$\sum_{i=0}^{n}\binom{n+1}{i}(n+i+1)B_{n+i} = 0.$$

Seidel's and Gessel's proofs are along similar lines although Seidel used a somewhat clumsy terminology.

These proofs of (2.4) show first that $B_{2n+1} = 0$ for $n \geq 1$. This is usually done by observing that (1.2) is equivalent with $e^z \sum_n B_n \frac{z^n}{n!} = z + \sum_n B_n \frac{z^n}{n!}$ or $B(z) = \sum_n B_n \frac{z^n}{n!} = \frac{z}{e^z - 1}$, which implies that $B(z) + \frac{z}{2} = \frac{z}{2} \frac{e^{\frac{z}{2}} + e^{-\frac{z}{2}}}{e^{\frac{z}{2}} - e^{-\frac{z}{2}}}$ is an even function.

To simplify the final step it is convenient to consider the linear functional $V$ on the polynomials defined by $V(x^n) = B_n$ for $n \neq 1$ and $V(x) = \frac{1}{2}$. This has the effect that in place of (1.2) we get

$$V\left((1-x)^n\right) = V(x^n) \tag{2.8}$$

for all $n \in \mathbb{N}$. Therefore by linearity

$$V\left(f(1-x)\right) = V\left(f(x)\right) \tag{2.9}$$

for each polynomial $f(x)$.

Choosing $f(x) = (1-x)^{n+1}(-x)^n$ we get $V\left(x^{n+1}(x-1)^n\right) = V\left((1-x)^{n+1}(-x)^n\right)$, which is the same as (2.7).

## 3. The Seidel triangle for Genocchi numbers

Seidel [9] has given a "Treppen-Schema" for the computation of the Genocchi numbers $G_{2n}$, which he called "Bernoulli'sche Zähler". We use a slightly changed version as in [3] and [10] and define the "Seidel triangle" for the Genocchi numbers as an array of integers $\left(g_{i,j}\right)_{i,j \geq 1}$ such that $g_{1,1} = g_{2,1} = 1$, $g_{i,j} = 0$ if $j < 0$ or $j > \left\lceil \frac{i}{2} \right\rceil$ and

$$g_{2i+1,j} = g_{2i+1,j-1} + g_{2i,j} \tag{3.1}$$

for $j = 1, 2, \cdots, i+1$ and

$$g_{2i,j} = g_{2i,j+1} + g_{2i-1,j} \tag{3.2}$$

for $j = i, i-1, \cdots, 1$.

$$\begin{pmatrix} 1 & 0 & 0 & 0 & 0 & 0 & 0 & 0 \\ 1 & 0 & 0 & 0 & 0 & 0 & 0 & 0 \\ 1 & 1 & 0 & 0 & 0 & 0 & 0 & 0 \\ 2 & 1 & 0 & 0 & 0 & 0 & 0 & 0 \\ 2 & 3 & 3 & 0 & 0 & 0 & 0 & 0 \\ 8 & 6 & 3 & 0 & 0 & 0 & 0 & 0 \\ 8 & 14 & 17 & 17 & 0 & 0 & 0 & 0 \\ 56 & 48 & 34 & 17 & 0 & 0 & 0 & 0 \end{pmatrix}$$

To compute this triangle note that the odd rows (3.1) are computed from left to right and the even rows (3.2) from right to left.
From (3.2) we see that

$$g_{2i,j} = \sum_{\ell \geq j} g_{2i-1,\ell} \tag{3.3}$$

and from (3.1) we have

$$g_{2i+1,j} = \sum_{\ell \leq j} g_{2i,\ell}. \tag{3.4}$$

We show now that the Seidel triangle is also closely related to the Fibonacci polynomials.

**Theorem 3.1**
*For $k = 1, 2, \cdots, n$*

$$g_{2n,k} = (-1)^n L(s^{n+1-k} F_{2k-1}) \tag{3.5}$$

*and for $k = 1, 2, \cdots, n+1$*

$$g_{2n+1,k} = (-1)^n L(s^{n+1-k} F_{2k}). \tag{3.6}$$

**Proof**

This is easily verified for $g_{1,1} = 1$ and $g_{2,1} = 1$.
The general case follows from
$g_{2n,k} = (-1)^n L(s^{n+1-k} F_{2k-1}) = (-1)^n L(s^{n-k} F_{2k+1}) - (-1)^n L(s^{n-k} F_{2k}) = g_{2n,k+1} + g_{2n-1,k}$
and
$g_{2n+1,k} = (-1)^n L(s^{n+1-k} F_{2k}) = (-1)^n L(s^{n-k+2} F_{2k-2}) + (-1)^n L(s^{n-k+1} F_{2k-1}) = g_{2n+1,k-1} + g_{2n,k}.$

As a special case we get that the "median Genocchi numbers" $H_{2n+1} := g_{2n+1,1}$ are given by

$$H_{2n+1} = (-1)^n L(s^n). \tag{3.7}$$

As another application we have

$$\sum_{j=0}^{k} \binom{k}{j} g_{n+j} = (-1)^{n+k-1} L\left( \sum_{j=0}^{k} (-1)^{k-j} \binom{k}{j} F_{n+j}(s) \right).$$

Observing that

$$\sum_{j=0}^{k}(-1)^{k-j}\binom{k}{j}F_{n+j}(s) = \frac{\alpha^n(\alpha-1)^k - \beta^n(\beta-1)^k}{\alpha-\beta} = s^k F_{n-k}(s)$$

we get

$$\sum_{j=0}^{k}\binom{k}{j}g_{n+j} = (-1)^{n+k-1}L\left(s^k F_{n-k}(s)\right). \tag{3.8}$$

**Corollary 3.2 (Seidel identity, L. Seidel [9], XIII) )**

$$\sum_{j=0}^{n}\binom{n}{j}g_{n+j} = \sum_{k=0}^{n}\binom{n}{2k}(-1)^k G_{2n-2k} = 0. \tag{3.9}$$

This is a counterpart to the v. Ettingshausen-Seidel-Kaneko identity (2.4).

**Corollary 3.3**

$$H_{2n+1} = (-1)^n \sum_{k=0}^{n+1}\binom{n+1}{k}g_{n+k} = \sum_{k=0}^{\frac{n}{2}}(-1)^k\binom{n+1}{2k+1}G_{2n-2k}. \tag{3.10}$$

By (3.8) and (3.7) and using $F_{-1}(s) = \frac{1}{s}$ we get

$$\sum_{k=0}^{n+1}\binom{n+1}{k}g_{n+k} = L(s^{n+1}F_{-1}(s)) = L(s^n) = (-1)^n H_{2n+1}.$$

These considerations also cast new light on [4], Theorem 1:

Apply $L$ to the generating function

$$\sum_{n\geq 0}F_{n+1}(s)z^n = \frac{1}{1-z-sz^2} = \frac{1}{1-z}\sum_{n\geq 0}s^n\left(\frac{z^2}{1-z}\right)^n.$$

This gives

$$1 + \sum_{n\geq 1}(-1)^{n-1}G_{2n}z^{2n-1} = \frac{1}{1-z}\sum_{n\geq 0}(-1)^n H_{2n+1}\left(\frac{z^2}{1-z}\right)^n. \tag{3.11}$$

which gives [4], Theorem 1, by multiplying with $z^2$ and $z \to -z$.

## 4. q-analogues

Next we show that the Seidel generation of the $q$–Genocchi numbers introduced in [10] is intimately related to the (Carlitz-) $q$–Fibonacci polynomials (cf. e.g. [1]).

Let

$$F_n(s,q) = \sum_{k=0}^{\frac{n-1}{2}} q^{2\binom{k}{2}} \begin{bmatrix} n-1-k \\ k \end{bmatrix} s^k. \tag{4.1}$$

Recall that

$$F_n(s,q) = F_{n-1}(s,q) + q^{n-3} s F_{n-2}(s,q). \tag{4.2}$$

Define a linear functional $L$ by

$$L\left(F_{2n+1}\left(s,\frac{1}{q}\right)\right) = [n=0]. \tag{4.3}$$

Now define a $q$–Seidel triangle $g_{n,k}(q)$ by

$$g_{2n,k}(q) = (-1)^n q^{2\binom{k-1}{2}} L\left(s^{n+1-k} F_{2k-1}\left(s,\frac{1}{q}\right)\right) \tag{4.4}$$

for $1 \leq k \leq n$
and

$$g_{2n+1,k}(q) = (-1)^n q^{(k-1)^2} L\left(s^{n+1-k} F_{2k}\left(s,\frac{1}{q}\right)\right) \tag{4.5}$$

for $1 \leq k \leq n+1$.
All other values should be $0$.

Then $g_{1,1}(q) = g_{2,1}(q) = 1$

and

$$g_{2n,k}(q) = g_{2n,k+1} + q^{k-1} g_{2n-1,k} \tag{4.6}$$

and

$$g_{2n+1,k}(q) = q^{k-1} g_{2n,k}(q) + g_{2n+1,k-1}(q). \tag{4.7}$$

This is precisely the definition given in [10].

The proof follows from

$$g_{2n,k}(q) = (-1)^n q^{2\binom{k-1}{2}} L\left(s^{n+1-k} F_{2k-1}\left(s,\frac{1}{q}\right)\right) = (-1)^n q^{2\binom{k-1}{2}} L\left(q^{2(k-1)}\left(s^{n-k} F_{2k+1}\left(s,\frac{1}{q}\right) - s^{n-k} F_{2k}\left(s,\frac{1}{q}\right)\right)\right)$$

$$= (-1)^n q^{2\binom{k}{2}} L\left(s^{n-k} F_{2k+1}\left(s,\frac{1}{q}\right)\right) + (-1)^{n-1} q^{(k-1)^2+(k-1)} L\left(s^{n-k} F_{2k}\left(s,\frac{1}{q}\right)\right) = g_{2n,k+1} + q^{k-1} g_{2n-1,k}$$

and

$$g_{2n+1,k}(q) = (-1)^n q^{(k-1)^2} L\left(s^{n+1-k} F_{2k}\left(s,\frac{1}{q}\right)\right)$$

$$= (-1)^n q^{(k-1)^2} L\left(s^{n+1-k} F_{2k-1}\left(s,\frac{1}{q}\right)\right) + (-1)^n q^{(k-1)^2 - 2k + 3} L\left(s^{n+2-k} F_{2k-2}\left(s,\frac{1}{q}\right)\right)$$

$$= q^{k-1} g_{2n,k}(q) + g_{2n+1,k-1}(q).$$

The first terms are

$$\begin{pmatrix} 1 & 0 & 0 & 0 \\ 1 & 0 & 0 & 0 \\ 1 & 1 & 0 & 0 \\ 1+q & q & 0 & 0 \\ 1+q & 1+q+q^2 & 1+q+q^2 & 0 \\ (1+q)^2(1+q^2) & q(1+q)(1+q+q^2) & q^2(1+q+q^2) & 0 \end{pmatrix}$$

Observe that from (4.7) we have

$$g_{2n+1,k}(q) = \sum_{\ell \geq 0} q^{k-1-\ell} g_{2n,k-\ell}(q) \tag{4.8}$$

and from (4.6)

$$g_{2n,k}(q) = \sum_{\ell \geq 0} q^{k-1+\ell} g_{2n-1,k+\ell}(q). \tag{4.9}$$

Now we follow [10] and define

$$G_{2n}(q) = g_{2n-1,n}(q) \tag{4.10}$$

and

$$H_{2n-1}(q) = q^{n-2} g_{2n-1,1}(q). \tag{4.11}$$

Then we get

$$L\left(sF_{2n}\left(s,\frac{1}{q}\right)\right) = (-1)^n g_{2n+1,n}(q) q^{-(n-1)^2} = (-1)^n q^{-(n-1)^2} G_{2n+2}(q), \tag{4.12}$$

$$L\left(F_{2n}\left(s,\frac{1}{q}\right)\right) = (-1)^{n-1} q^{-(n-1)^2} G_{2n}(q) \tag{4.13}$$

and

$$L(s^n) = (-1)^n g_{2n+1,1}(q) = (-1)^n \frac{H_{2n+1}(q)}{q^{n-1}}. \tag{4.14}$$

Next we prove

$$\sum_{k=0}^{n}(-1)^k \begin{bmatrix} n \\ k \end{bmatrix} q^{\binom{k}{2}} F_{2n-k}(s,q) = 0. \tag{4.15}$$

By changing $q \to \frac{1}{q}$ we get

$$\sum_{k=0}^{n}(-1)^k \begin{bmatrix} n \\ k \end{bmatrix} q^{\binom{k+1}{2}-kn} F_{2n-k}\left(s,\frac{1}{q}\right) = 0. \tag{4.16}$$

If we apply the linear functional $L$ and observe (4.13) we get

$$(-1)^{2n-1} \begin{bmatrix} n \\ 2n-1 \end{bmatrix} q^{\binom{2n}{2}-(2n-1)n} F_1(s,\frac{1}{q}) + \sum_{k=0}^{n}(-1)^{n-k-1} \begin{bmatrix} n \\ 2k \end{bmatrix} q^{\binom{2k+1}{2}-2kn-(n-k-1)^2} G_{2n-2k}(q) = 0.$$

Therefore we get
**Theorem 4.1 (q-Seidel identity)**

$$\sum_{k=0}^{n}(-1)^k \begin{bmatrix} n \\ 2k \end{bmatrix} q^{2\binom{k}{2}} G_{2n-2k}(q) = [n=1]. \tag{4.17}$$

This formula can be used to compute the polynomials $G_{2n}(q)$.

In order to prove (4.15) we show more generally

$$\sum_{k=0}^{n}(-1)^k q^{\binom{k}{2}} \begin{bmatrix} n \\ k \end{bmatrix} F_{2n+m-k}(s,q) = q^{2\binom{n}{2}+(m-1)n} s^n F_m(s,q). \tag{4.18}$$

For $n = 0$ this is trivial.
For $n = 1$ (4.18) reduces to the recursion

$$F_{m+2}(s,q) - F_{m+1}(s,q) = q^{m-1} s F_m(s), \tag{4.19}$$

which holds for $m \in \mathbb{Z}$.

Suppose that (4.18) is already known for $i < n$ and all $m \in \mathbb{Z}$.
Then

$$\sum_k (-1)^k q^{\binom{k}{2}} \begin{bmatrix} n \\ k \end{bmatrix} F_{2n+m-k}(s,q) = \sum_k (-1)^k q^{\binom{k}{2}} q^k \begin{bmatrix} n-1 \\ k \end{bmatrix} F_{2n+m-k}(s,q) + \sum_k (-1)^k q^{\binom{k}{2}} \begin{bmatrix} n-1 \\ k-1 \end{bmatrix} F_{2n+m-k}(s,q)$$

$$= \sum_k (-1)^{k-1} q^{\binom{k}{2}} \begin{bmatrix} n-1 \\ k-1 \end{bmatrix} F_{2n+m-k+1}(s,q) + \sum_k (-1)^k q^{\binom{k}{2}} \begin{bmatrix} n-1 \\ k-1 \end{bmatrix} F_{2n+m-k}(s,q)$$

$$= \sum_k (-1)^{k-1} q^{\binom{k}{2}} \begin{bmatrix} n-1 \\ k-1 \end{bmatrix} (F_{2n+m-k+1}(s,q) - F_{2n+m-k}(s,q))$$

$$= \sum_k (-1)^k q^{\binom{k+1}{2}} \begin{bmatrix} n-1 \\ k \end{bmatrix} q^{2n+m-k-3} s F_{2n+m-k-2}(s,q) = q^{2n+m-3} s \sum_k (-1)^k q^{\binom{k}{2}} \begin{bmatrix} n-1 \\ k \end{bmatrix} F_{2(n-1)+m-k}(s,q)$$

$$= q^{2n+m-3} q^{2\binom{n-1}{2}+(m-1)(n-1)} s s^{n-1} F_m(s,q) = q^{2\binom{n}{2}+(m-1)n} s^n F_m(s,q).$$

A $q$-analogue of (3.10) is

$$H_{2n+1}(q) = \sum_{k=0}^{n} (-1)^k \begin{bmatrix} n+1 \\ 2k+1 \end{bmatrix} q^{k^2-k+n-2} G_{2n-2k}(q) \qquad (4.20)$$

for $n \geq 2$.

If we choose $m = -1$ in (4.18) and replace $n$ by $n+1$ we get

$$\sum_{k=0}^{n} (-1)^k q^{\binom{k}{2}} \begin{bmatrix} n+1 \\ k \end{bmatrix} F_{2(n+1)-1-k}(s,q) = q^{2\binom{n+1}{2}+(-1-1)(n+1)} s^{n+1} F_{-1}(s,q)$$

$$= q^{n^2-n-2} s^{n+1} \frac{q^2}{s} = q^{2\binom{n}{2}} s^n.$$

If we now change $q \to \frac{1}{q}$ we have

$$s^n = q^{2\binom{n}{2}} \sum_{k=0}^{n} (-1)^k q^{\binom{k}{2}-kn} \begin{bmatrix} n+1 \\ k \end{bmatrix} F_{2n+1-k}\left(s, \frac{1}{q}\right).$$

By applying the linear functional $L$ we obtain the desired result

$$H_{2n+1}(q) = (-1)^n q^{n-1} L(s^n) = (-1)^n q^{n-1} q^{2\binom{n}{2}} \sum_{k=0}^{n} (-1)^k q^{\binom{k}{2}-kn} \begin{bmatrix} n+1 \\ k \end{bmatrix} L\left(F_{2n+1-k}\left(s, \frac{1}{q}\right)\right)$$

$$= (-1)^{n-1} q^{n-1} q^{2\binom{n}{2}} \sum_{k=0}^{n} q^{\binom{2k+1}{2}-(2k+1)n} \begin{bmatrix} n+1 \\ 2k+1 \end{bmatrix} L\left(F_{2n+1-2k-1}\left(s, \frac{1}{q}\right)\right)$$

$$= \sum_{k=0}^{n} (-1)^k q^{k^2-k+n-2} \begin{bmatrix} n+1 \\ 2k+1 \end{bmatrix} G_{2n-2k}(q).$$

There is also a $q$-analogue of (3.11).

In the generating function

$$\sum_n F_{n+1}(s,q)z^n = \sum_n z^n \sum_{k=0}^{n} q^{2\binom{k}{2}} \begin{bmatrix} n-k \\ k \end{bmatrix} s^k = \sum_k q^{2\binom{k}{2}} s^k z^{2k} \sum_j \begin{bmatrix} k+j \\ k \end{bmatrix} z^j = \sum_k q^{2\binom{k}{2}} \frac{s^k z^{2k}}{(1-z)\cdots(1-q^k z)}$$

change $q$ to $\frac{1}{q}$. We obtain

$$\sum_n F_{n+1}\left(s,\frac{1}{q}\right)z^n = \sum_k (-1)^{k+1} q^{-2\binom{k}{2}+\binom{k+1}{2}} \frac{s^k z^{2k}}{(z-1)\cdots(z-q^k)}.$$

By applying $L$ we get

$$1 + \sum_{n \geq 1}(-1)^{n-1} q^{-(n-1)^2} G_{2n}(q) z^{2n-1} = -\sum_k q^{1-\binom{k}{2}} \frac{z^{2k}}{(z-1)\cdots(z-q^k)} H_{2k+1}(q)$$

$$= \sum_k q^{k-\binom{k}{2}} (-1)^k \frac{z^{2k}}{(1-z)\cdots(q^k-z)} g_{2k+1,1}(q). \tag{4.21}$$

**Remark**

For $q = 1$ some of these results have been obtained by using Seidel matrices instead of the above Seidel triangle.

The Seidel matrix $(a_{n,k})$ for a given sequence $(c_n)$ is defined by $a_{n,0} = c_n$ and

$a_{n,k} = a_{n,k-1} + a_{n+1,k-1}$ for $k \geq 1$. Then $a_{n,k} = \sum_{i=0}^{k} \binom{k}{i} a_{n+i,0}$.

To obtain a useful $q$-analogue define instead

$$a_{n,k} = q^{n-1}\left(a_{n,k-1} + a_{n+1,k-1}\right). \tag{4.22}$$

Then

$$a_{n,k} = q^{k(n-1)} \sum_{j=0}^{k} q^{\binom{j}{2}} \begin{bmatrix} k \\ j \end{bmatrix} a_{n+j,0}. \tag{4.23}$$

This holds for $k = 1$.
By induction we get

$$a_{n,k+1} = q^{n-1}\left(a_{n,k} + a_{n+1,k}\right) = q^{n-1}\left(q^{k(n-1)} \sum_{j=0}^{k} q^{\binom{j}{2}} \begin{bmatrix} k \\ j \end{bmatrix} a_{n+j,0} + q^{kn} \sum_{j=0}^{k} q^{\binom{j-1}{2}} \begin{bmatrix} k \\ j-1 \end{bmatrix} a_{n+j,0}\right)$$

$$= q^{n-1+kn-k} \sum_j q^{\binom{j}{2}} a_{n+j,0} \left(\begin{bmatrix} k \\ j \end{bmatrix} + q^{k-j+1} \begin{bmatrix} k \\ j-1 \end{bmatrix}\right) = q^{(k+1)(n-1)} \sum_j q^{\binom{j}{2}} \begin{bmatrix} k+1 \\ j \end{bmatrix} a_{n+j,0}.$$

By choosing

$$c_n = a_{n,0} = L\left((-1)^{n-1} F_n\left(s,\frac{1}{q}\right)\right) \tag{4.24}$$

we get

$$a_{n,k} = L\left((-1)^{n-k-1} q^{\binom{k+1}{2}} s^k F_{n-k}\left(s, \frac{1}{q}\right)\right). \tag{4.25}$$

I want also mention a $q$ – analogue of the v. Ettingshausen- Seidel-Kaneko identity.
If we apply the linear functional $M$ defined by $M(F_{2n+2}(s,q)) = [n=0]$ to (4.15) we get

$$\sum_{k=0}^{n+1} (-1)^k \begin{bmatrix} n+1 \\ k \end{bmatrix} q^{\binom{k}{2}} M(F_{2n+2-k}(s,q)) = 0. \tag{4.26}$$

This can be used to compute the sequence $(M(F_{2k+1}(s,q)))$ which begins with

$$\left\{1, \frac{q}{1+q}, -\frac{q^4}{(1+q)(1+q+q^2)}, \frac{q^7}{(1+q)(1+q+q^2)}, \right.$$
$$\left. -\frac{q^{10}(1+q+2q^2+2q^3+q^4+q^5+q^6)}{(1+q)(1+q+q^2)(1+q+q^2+q^3+q^4)}, \frac{q^{13}(1+q+3q^2+3q^3+2q^4+2q^5+q^6+q^7+q^8)}{(1+q)(1+q+q^2)^2}\right\}$$

## 5. Some identities

Each identity for the $q$ – Fibonacci polynomials gives an identity for the entries of the
$q$ – Seidel triangle $(g_{i,j}(q))$.
I shall give some examples.

1) From the definition of the $q$ – Fibonacci polynomials we get

$$F_n\left(s, \frac{1}{q}\right) = \sum_{k=0}^{\frac{n-1}{2}} q^{k^2+2k-nk} \begin{bmatrix} n-1-k \\ k \end{bmatrix} s^k.$$

Applying the linear functional $L$ gives

$$\sum_{k=0}^{n} (-1)^k q^{k^2+k-2nk} \begin{bmatrix} 2n-k \\ k \end{bmatrix} g_{2k+1,1}(q) = 0 \tag{5.1}$$

and

$$q^{(n-1)^2} \sum_{k=0}^{n-1} (-1)^{n-k-1} q^{k^2+2k-2nk} \begin{bmatrix} 2n-1-k \\ k \end{bmatrix} g_{2k+1,1}(q) = G_{2n}(q). \tag{5.2}$$

These identities are $q$ – analogues of [3], Corollary 1.

2) It is easily verified (cf. [2]) that

$$F_{m+2n}(s,q) = \sum_{k=0}^{n} \begin{bmatrix} n \\ k \end{bmatrix} q^{k(m+n-2)} s^k F_{m+n-k}(s,q). \tag{5.3}$$

Therefore

$$F_{2n+1}(s,q) = \sum_{j=0}^{\frac{n}{2}} \begin{bmatrix} n \\ n-2j \end{bmatrix} q^{(n-2j)(n-1)} s^{n-2j} F_{2j+1}(s,q) + \sum_{j=0}^{\frac{n+1}{2}} \begin{bmatrix} n \\ n-2j+1 \end{bmatrix} q^{(n-2j+1)(n-1)} s^{n-2j+1} F_{2j}(s,q)$$

By changing $q \to \dfrac{1}{q}$ we get

$$F_{2n+1}\left(s,\frac{1}{q}\right) = \sum_{j} \begin{bmatrix} n \\ 2j \end{bmatrix} q^{2\binom{2j}{2}-2\binom{n}{2}} s^{n-2j} F_{2j+1}\left(s,\frac{1}{q}\right) + \sum_{j} \begin{bmatrix} n \\ 2j-1 \end{bmatrix} q^{2\binom{2j-1}{2}-2\binom{n}{2}} s^{n-2j+1} F_{2j}\left(s,\frac{1}{q}\right)$$

If we apply the linear functional $L$ we get

$$0 = \sum_{j} (-1)^j \begin{bmatrix} n \\ 2j \end{bmatrix} q^{2\binom{2j}{2}-2\binom{j}{2}} g_{2n-2j,j+1}(q) + \sum_{j} (-1)^j \begin{bmatrix} n \\ 2j-1 \end{bmatrix} q^{2\binom{2j-1}{2}-(j-1)^2} g_{2n-2j+1,j}(q)$$

or

$$\sum_{j}(-1)^j \begin{bmatrix} n \\ 2j \end{bmatrix} q^{3j^2-j} g_{2n-2j,j+1}(q) = \sum_{j} (-1)^{j-1} \begin{bmatrix} n \\ 2j-1 \end{bmatrix} q^{3j^2-4j+1} g_{2n-2j+1,j}(q). \qquad (5.4)$$

3) Now consider the following identities which in the special case $m = -1$ have been proved in [8]:

$$\sum_{k=0}^{n} (-s)^k \begin{bmatrix} 2k+m \\ k \end{bmatrix} q^{-\binom{k+m+2}{2}} = \sum_{k=0}^{n} (-s)^{n-k} \begin{bmatrix} 2n+m+2 \\ n-k \end{bmatrix} q^{-\binom{n+2+m-k}{2}} F_{2k+2}(s,q) \qquad (5.5)$$

and

$$\sum_{k=0}^{n} (-s)^k \begin{bmatrix} 2k+m \\ k \end{bmatrix} q^{-\binom{k+m+2}{2}} = \sum_{k=0}^{n} (-s)^{n-k} \begin{bmatrix} 2n+m+1 \\ n-k \end{bmatrix} q^{-\binom{n+2+m-k}{2}} F_{2k+1}(s,q). \qquad (5.6)$$

Comparing coefficients equation (5.5) is equivalent with

$$(-1)^k \begin{bmatrix} 2k+m \\ k \end{bmatrix} q^{-\binom{k+m+2}{2}} = \sum_{i+j=k} (-1)^j \begin{bmatrix} 2n+m+2 \\ j \end{bmatrix} q^{-\binom{j+2+m}{2}} \begin{bmatrix} 2n-2j+1-i \\ i \end{bmatrix} q^{i^2-i}$$

Observing that $\begin{bmatrix} -r \\ k \end{bmatrix} = (-1)^k q^{-\left(kr+\binom{k}{2}\right)} \begin{bmatrix} r+k-1 \\ k \end{bmatrix}$ we get

$$\sum_{i+j=k} (-1)^j \begin{bmatrix} 2n+m+2 \\ j \end{bmatrix} q^{-\binom{j+2+m}{2}} \begin{bmatrix} 2n-2j+1-i \\ i \end{bmatrix} q^{2\binom{i}{2}}$$

$$= (-1)^k \sum_{j} \begin{bmatrix} 2n+m+2 \\ j \end{bmatrix} q^{-\binom{j+2+m}{2}} q^{2\binom{k-j}{2}} \begin{bmatrix} 2k-2n-2 \\ k-j \end{bmatrix} q^{(k-j)(2n-j-k+1)-\binom{k-j}{2}}.$$

This can be simplified to $\begin{bmatrix} 2k+m \\ k \end{bmatrix} = \sum_j q^{(k-j)(2n+m+2-j)} \begin{bmatrix} 2n+m+2 \\ j \end{bmatrix} \begin{bmatrix} 2k-2n-2 \\ k-j \end{bmatrix}$ which is true by the $q$-Vandermonde theorem.

In the same way we can prove (5.6).

By changing $q \to \frac{1}{q}$ we get

$$\sum_{k=0}^n (-s)^k \begin{bmatrix} 2k+m \\ k \end{bmatrix} \frac{1}{q^{\frac{k(k-3)}{2}}} = \frac{1}{q^{\binom{n+1}{2}}} \sum_{k=0}^n (-s)^{n-k} \begin{bmatrix} 2n+m+2 \\ n-k \end{bmatrix} q^{\frac{k(3k+1)}{2}} q^{-kn} F_{2k+2}\left(s, \frac{1}{q}\right) \quad (5.7)$$

and

$$\sum_{k=0}^n (-s)^k \begin{bmatrix} 2k+m \\ k \end{bmatrix} \frac{1}{q^{\frac{k(k-3)}{2}}} = \frac{1}{q^{\binom{n}{2}}} \sum_{k=0}^n (-s)^{n-k} \begin{bmatrix} 2n+m+1 \\ n-k \end{bmatrix} q^{\frac{k(3k-1)}{2}} q^{-kn} F_{2k+1}\left(s, \frac{1}{q}\right). \quad (5.8)$$

Applying $L$ we get the identity

$$\sum_{k=0}^n \begin{bmatrix} 2k+m \\ k \end{bmatrix} \frac{1}{q^{\frac{k(k-3)}{2}}} g_{2k+1,1}(q) = \frac{1}{q^{\binom{n+1}{2}}} \sum_{k=0}^n (-1)^k \begin{bmatrix} 2n+m+2 \\ n-k \end{bmatrix} q^{\binom{k+1}{2}-kn} g_{2n+1,k+1}(q)$$
$$= \frac{1}{q^{\binom{n}{2}}} \sum_{k=0}^n (-1)^k \begin{bmatrix} 2n+m+1 \\ n-k \end{bmatrix} q^{\binom{k+1}{2}-kn} g_{2n,k+1}(q). \quad (5.9)$$

For $n=1$ this reduces to
$$1 + q[m+2] = \frac{[m+4]-1}{q} = [m+3]$$

and for $n=2$ we get

$$1 + q[m+2] + q\begin{bmatrix} m+4 \\ 2 \end{bmatrix}[2] = \begin{bmatrix} m+6 \\ 2 \end{bmatrix}\frac{[2]}{q^3} - \begin{bmatrix} m+6 \\ 1 \end{bmatrix}\frac{[3]}{q^4} + \frac{[3]}{q^4} = \begin{bmatrix} m+5 \\ 2 \end{bmatrix}\frac{[2]}{q} - \begin{bmatrix} m+5 \\ 1 \end{bmatrix}\frac{1}{q}.$$